# 拟连续 domain 上的广义理想收敛

王武[1]，谭彬[2]，张舜[3]

（1. 天津理工大学中环信息学院 基础课部，天津 300380；2.天津理工大学 理学院，天津 300384；3. 天津仁爱学院 数理教学部，天津 300163）

**摘　要**：在定向完备偏序集中引入了广义理想下极限和广义理想终下界极限的概念，并研究了其与 Scott 拓扑和 Lawson 拓扑的关系．主要结果有：(1)在定向完备偏序集上，广义理想下极限拓扑与 Scott 拓扑一致；(2) 广义理想下极限收敛是拓扑的当且仅当定向完备偏序集是拟连续 domain；(3)在拟连续 domain 中，广义理想终下界极限拓扑与 Lawson 拓扑一致，并给出了定向完备偏序集为连续 domain 的一个充分条件．

**关键词**： 拟连续偏序集；广义理想下极限；广义理想终下界极限；Scott拓扑；Lawson拓扑

**中图分类号**：O153.1；O189.1　　**文献标志码**：A　　**文章编号**：

## Generalized ideal convergence on quasi-continuous domains

WANG Wu[1], TAN Bin[2], ZHANG Shun[3]

(1.Basic Course Department of Zhonghuan Information College Tianjin University of Technology，tianjin 300380，China；
2. College Science of Tianjin University of Technology，tianjin 300384，China;
3. Mathematics Teaching Department of Tianjin Ren'ai College，tianjin 301636，China)

**Abstract:** In this paper，the concepts of generalized ideal inf-limit and generalized ideal final lower bound limit are introduced in the directed complete poset, and their relations with Scott topology and Lawson topology are studied. The main results are as follows: (1) On directed complete posets, generalized ideal inf-limit topology is consistent with Scott topology; (2) Generalized ideal inf-limitis convergence is topological if and only if directed complete posets are quasi-continuous domains; (3) In quasi-continuous domain, generalized ideal final lower bound limit topology is consistent with Lawson topology, and a sufficient condition of continuous domain is given.

**Key words:** quasi-continuous domain；generalized ideal inf-limit；generalized ideal final lower bound limit；Scott topology；Lawson topology

　　偏序集理论作为重要的数学分支之一，旨在为计算机高级程序设计语言提供数学模型，受到了计算机科学和数学领域诸多学者的关注，并得到了很多有价值的结论和模型[1-6]．随着计算机理论的发展，偏序集理论不断向信息科学、逻辑学、分析学及各种应用学科渗透[7-8]．Domain 和拟连续 Domain 作为特殊的偏序结构，有很多良好的性质和特征，并且在理论计算机中的应用更加广泛[9]，如李高林等给出了拟连续 domain 的拟基的概念，并研究了拟基的一些性质[10]；赵浩然等研究了拟连续 domain 的函数空间[11]，并说明了由拟连续 domain 作为对象，Scott 连续函数作为映射的范畴不是 Cartesian 闭范畴；徐晓泉[12]等研究了拟连续 domain 上 Lawson 拓扑的紧性．网的收敛是拓扑空间中的重要工具，可以完全刻画拓扑空间。本文利用定向集的理想定义了网的广义理想下极限和广义理想终下界极限，并研究了其与 Scott 拓扑和 Lawson 拓扑的关系．主要结果有：(1)在定向完备偏序集上，广义理想下极限拓扑与 Scott 拓扑一致；(2) 广义理想下极限收敛是拓扑的当且仅当定向完备偏序集是拟连续 domain；(3)在拟连续 domain 中，广义理想终下界极限拓扑与 Lawson 拓扑一致，并给出了定向完备偏序集为连续 domain 的等价刻画．本文结果有助于 domain 理论的进一步研究．

## 1　预备知识

首先，介绍偏序集中的一些基本概念[12]. 设 $L$ 是偏序集，$D \subseteq L$，如果 $D$ 中任意两个元在 $D$ 中有上界，则称 $D$ 为定向集. 如果 $L$ 的每个定向子集 $D$ 都有上确界(记为 $\sup D$)，则称 $L$ 为定向完备偏序集(简称 $dcpo$). 任给 $A \subseteq L$，记 $\uparrow A = \{x \in L : \exists a \in A, a \leq x\}$，$\downarrow A = \{x \in L : \exists a \in A, x \leq a\}$；若 $A$ 为单点集 $\{a\}$，则有记号 $\uparrow A = \uparrow a$，$\downarrow A = \downarrow a$. 设 $L$ 是定向完备偏序集，$x, y \in L$，$D \subseteq L$ 是定向集，如果 $y \leq \sup D$ 蕴含 $D \cap \uparrow x \neq \varnothing$，则 $x \ll y$，记 $\Downarrow x = \{y : y \ll x\}$. 如果任意 $x \in L$，集合 $\Downarrow x$ 是定向的并且 $x = \vee \Downarrow x$，则称 $L$ 为连续的，称连续的定向完备偏序集为 domain 或连续 domain.

设 $L$ 是定向完备偏序集，$G, H$ 是 $L$ 的两个非空子集. 如果 $H \subseteq \uparrow G$，则 $G \leq H$，称这种序关系为 Symth 序，Symth 序是一种预序关系. 如果对任意的定向集 $D \subseteq L$，$\sup D \in \uparrow H$ 蕴含存在 $d \in D$ 使得 $d \in \uparrow G$，则称 $G$ 逼近 $H$，记为 $G \ll H$. 特别的，$\{y\} \ll H$ 简记为 $y \ll H$，$G \ll \{x\}$ 简记为 $G \ll x$. 显然，$G \ll H$ 蕴含 $\forall h \in H$，$G \ll h$. 易知上述定义的逼近关系有如下性质：(1) $G \ll H$ 蕴含 $G \leq H$；(2) $G \leq E \ll F \leq H$ 蕴含 $G \ll H$.

令 $\mathcal{P}^w(L)$ 表示 $L$ 的所有非空有限子集构成的集族. 设非空集族 $\mathcal{F} \subseteq \mathcal{P}^w(L)$，如果任给 $E, F \in \mathcal{F}$ 存在有限集 $H \in \mathcal{F}$ 满足 $\uparrow H \subseteq \uparrow E \cap \uparrow F$，则称 $\mathcal{F}$ 为定向集族. 设 $L$ 是定向完备偏序集，如果任意 $x \in L$，集族 $fin(x) = \{F \in \mathcal{P}^w(L) : F \ll x\}$ 是定向集族且 $\uparrow x = \bigcap_{F \in fin(x)} \uparrow F$，则称 $L$ 是拟连续 domain. 拟连续 domain 的等价刻画是：存在定向集族 $\mathcal{F} \subseteq fin(x)$ 使得 $\uparrow x = \bigcap_{F \in \mathcal{F}} \uparrow F$.

**定义 1**[13]　设 $L$ 是偏序集.

(1) 任给子集 $U \subseteq L$，如果 $U = \uparrow U$ 且对任意的定向集 $D \subseteq L$，$\sup D \in U$ 蕴含 $D \cap U \neq \varnothing$，则称 $U$ 为 Scott 开集. 所有的 Scott 开集构成一个拓扑，称为 Scott 拓扑，记为 $\sigma(L)$；

(2) 以主滤子的补 $L \setminus \uparrow x$ 为子基生成的拓扑称为下拓扑，记为 $\omega(L)$；

(3) 称由 $\sigma(L)$ 和 $\omega(L)$ 共同加细的拓扑 $\sigma(L) \vee \omega(L)$ 为 Lawson 拓扑，记为 $\lambda(L)$.

设 $L$ 是定向完备偏序集，$x \in L$，$D \subseteq L$ 是定向集，如果 $x \leq \sup D$ 蕴含 $x \in cl_\sigma(\downarrow x \cap \downarrow D)$，则称 $L$ 是交连续的. $L$ 是交连续的当且仅当任意 $x \in L$，$U \in \sigma(L)$，$\uparrow(U \cap \downarrow x) \in \sigma(L)$. 定向完备偏序集是连续的当且仅当它既是交连续的又是拟连续的[13].

**命题 1**[13]　设 $L$ 是拟连续 domain，记 $\Uparrow F = \{x : F \ll x\}$，则

(1) 任意 $x \in L$，$H$ 为 $L$ 的有限子集。如果 $H \ll x$，则存在有限集 $F \subseteq L$ 使得 $H \ll F \ll x$；

(2) $\mathrm{int}_{\sigma(L)}(\uparrow F) = \Uparrow F$，其中 $\mathrm{int}_{\sigma(L)}(\uparrow F)$ 表示 $\uparrow F$ 在 Scott 拓扑 $\sigma(L)$ 中的内部；

(3) $\{\Uparrow F : F \in \mathcal{P}^w(L)\}$ 为 Scott 拓扑 $\sigma(L)$ 的基.

设 $M$ 是一个集合，如果 $M$ 的子集族 $\mathcal{I}$ 满足下列条件：(1) $A \in \mathcal{I}$，$B \subseteq A$ 蕴含 $B \in \mathcal{I}$；(2) $A, B \subseteq \mathcal{I}$ 蕴含 $A \cup B \in \mathcal{I}$，则称 $\mathcal{I}$ 是 $M$ 的一个理想. 如果 $M \notin \mathcal{I}$，则称 $\mathcal{I}$ 是 $M$ 的非平凡理想. 令 $J$ 是定向集，$M_j = \{j' \in J : j' \in \uparrow j\}$，则 $\mathcal{I}_0 = \{A \subseteq J : A \subseteq J \setminus M_j\}$ 是非平凡理想. 显然任意有限集 $E \subseteq M$，$e \in E$，$M_e \subseteq \{j' \in J : j' \in \uparrow E\} = M_E$，则 $A = J \setminus M_E \subseteq J \setminus M_e$，则 $J \setminus M_E \in \mathcal{I}_0$.

**定义 2**[14]　设 $L$ 是偏序集，$(x_j)_{j \in J}$ 是 $L$ 中的一个网，$\mathcal{I}$ 是 $J$ 的理想，$x \in L$. 如果存在定向集 $D \subseteq L$ 使得下列条件成立：(1) $x \leq \sup D$；(2)任意的 $d \in D$，$\{j \in J : x_j \not\geq d\} \in \mathcal{I}$. 则称网 $(x_j)_{j \in J}$ 是理想下极限收敛到 $x$. 也称 $x$ 是 $(x_j)_{j \in J}$ 的理想下极限，记为 $(x_j)_{j \in J} \to_{\mathcal{IS}} x$. 设 $\mathcal{I}$ 是 $J$ 的理想，令
$$\tau_L^{I-\lim-\inf} = \{U \subseteq L : (x_j)_{j \in J} \to_{\mathcal{IS}} x, \ x \in U, \ \text{则} \ \{j \in J : x_j \notin U\} \in \mathcal{I}\}.$$

**命题 2**[14]　设 $L$ 是偏序集，则 $\sigma_L = \tau_L^{I-\lim-\inf}$.

显然，如果 $L$ 是定向完备偏序集，$\sigma_L = \tau_L^{I-\lim-\inf}$.

**定义 3**[14]　设 $X$ 是拓扑空间，$(x_j)_{j \in J}$ 是 $L$ 中的一个网，$\mathcal{I}$ 是 $J$ 的理想，若 $x \in U$，$\{j \in J : x_j \notin U\} \in \mathcal{I}$，则称网 $(x_j)_{j \in J}$ 是理想收敛到 $x$ 的. 也称 $x$ 是 $(x_j)_{j \in J}$ 的理想极限，记为 $(x_j)_{j \in J} \to_{\mathcal{I}} x$.

易知，拓扑空间 $X$ 中的网 $(x_j)_{j \in J}$ 收敛到 $x$ 当且仅当 $(x_j)_{j \in J} \to_{\mathcal{I}_0} x$.

**命题 3**[15]　**(Rudin 引理)**　设 $L$ 是偏序集，$\mathcal{F} \subseteq \mathcal{P}^w(L)$ 是定向集族. 那么存在定向集 $D \subseteq \bigcup_{F \in \mathcal{F}} F$ 满足 $D \cap F \neq \varnothing$ 对所有 $F \in \mathcal{F}$ 成立.

Rudin 引理的一个直接推论是：设 $L$ 是偏序集，$\mathcal{F} \subseteq \mathcal{P}^w(L)$ 是定向集族，$U \subseteq L$ 是 Scott 开集。则 $\bigcap_{F \in \mathcal{F}} F \subseteq U$ 意味存在 $F \in \mathcal{F}$ 满足 $\uparrow F \subseteq U$。

## 2 广义理想下极限收敛与 Scott 拓扑

本节将连续 domain 的理想下极限收敛进行推广，在定向完备偏序集中引入广义理想下极限收敛的概念。并证明广义理想下极限收敛是拓扑的当且仅当定向完备偏序集是拟连续 domain。

**定义 4** 设 $L$ 是定向完备偏序集，$(x_j)_{j \in J}$ 是 $L$ 中的一个网，$\mathcal{I}$ 是 $J$ 的理想，$x \in L$。如果存在定向集族 $\mathcal{F} \subseteq \mathcal{P}^w(L)$ 使得下列条件成立：(1) $\bigcap_{F \in \mathcal{F}} \uparrow F \subseteq \uparrow x$；(2) 任意的 $F \in \mathcal{F}$，$\{j \in J : x_j \notin \uparrow F\} \in \mathcal{I}$。则称网 $(x_j)_{j \in J}$ 广义理想下极限收敛到 $x$，也称 $x$ 是 $(x_j)_{j \in J}$ 的广义理想下极限，记为 $(x_j)_{j \in J} \to_{\mathcal{GIS}} x$。

设 $\mathcal{I}$ 是 $J$ 的理想，令

$$\mathcal{GIS}(L) = \{((x_j)_{j \in J}, x) : (x_j)_{j \in J} \to_{\mathcal{GIS}} x\},$$

$$\tau_L^{g-I-\lim\text{-inf}} = \{U \subseteq L : ((x_j)_{j \in J}, x) \in \mathcal{GIS}(L),\ x \in U,\ \text{则}\ \{j \in J : x_j \notin U\} \in \mathcal{I}\}.$$

**命题 4** 设 $L$ 是定向完备偏序集，$(x_j)_{j \in J}$ 是 $L$ 中的一个网，$\mathcal{I}$ 是 $J$ 的理想，$x \in L$。则 $(x_j)_{j \in J} \to_{\mathcal{IS}} x$ 蕴含 $(x_j)_{j \in J} \to_{\mathcal{GIS}} x$。

**证明** 设 $(x_j)_{j \in J} \to_{\mathcal{IS}} x$，则存在定向集 $D \subseteq L$ 使得 $x \leq \sup D$ 且任意的 $d \in D$，$\{j \in J : x_j \not\geq d\} \in \mathcal{I}$。令 $\mathcal{F} = \{\{d\} : d \in D\}$，则 $\uparrow \sup D = \bigcap_{F \in \mathcal{F}} \uparrow F$。由于 $x \leq \sup D$，则 $\bigcap_{F \in \mathcal{F}} \uparrow F \subseteq \uparrow x$。并且任意的 $F \in \mathcal{F}$，$\{j \in J : x_j \notin \uparrow F\} = \{j \in J : x_j \not\geq d\} \in \mathcal{I}$ 成立，故 $(x_j)_{j \in J} \to_{\mathcal{GIS}} x$。

上述结果的逆命题不一定成立，从而说明广义理想下极限收敛是理想下极限收敛的推广。

**例 1** 设 $L = N \cup \{\infty, a\}$，这里 $N$ 表示自然数的集合。任意 $x, y \in L$，$x \leq y \Leftrightarrow y = \infty$ 或 $x, y \in N$，$x \leq y$。易知，$L$ 是非连续的拟连续 domain，满足任意的 $n \in N$，$\{a, n\} \ll a$ 且 $\uparrow a = \bigcap_{n \in N} \uparrow \{a, n\}$。令 $x_{2n} = n$，$x_{2n+1} = a$，则 $(x_n)_{n \in N}$ 是一个网。取 $N$ 的理想 $\mathcal{I}_0$，则显然 $\{m : x_m \notin \uparrow \{n, a\}\} \in \mathcal{I}_0$，令 $\mathcal{F} = \{\{n, a\} : n \in N\}$，则 $\bigcap_{F \in \mathcal{F}} \uparrow F = \{a, \infty\} \subseteq \uparrow a$，$\{m : x_m \notin \uparrow \{n, a\}\} \in \mathcal{I}_0$，即 $(x_n)_{n \in N} \to_{\mathcal{GIS}} a$。另一方面，如果 $x \in L$，则显然 $\{m : x_m \not\geq \uparrow n\} \notin \mathcal{I}_0$，从而不存在定向集 $D \subseteq L$ 使得 $x \leq \sup D$ 且任意的 $d \in D$，$\{j \in J : x_j \not\geq d\} \in \mathcal{I}_0$，即 $(x_n)_{n \in N} \to_{\mathcal{IS}} a$ 不成立。

**命题 5** 设 $L$ 是定向完备偏序集，$G \in \mathcal{P}^w(L)$，$x \in L$。如果任意的网 $(x_j)_{j \in J} \subseteq L$ 以及 $J$ 的非平凡理想 $\mathcal{I}$，$(x_j)_{j \in J} \to_{\mathcal{GIS}} x$ 蕴含 $\{j \in J : x_j \notin \uparrow G\} \in \mathcal{I}$，则 $G \ll x$。

**证明** 设定向集 $D \subseteq L$，$x \leq \sup D$，令 $(x_d)_{d \in D}$ 满足 $x_d = d$，则 $(x_d)_{d \in D} \to_{\mathcal{I}_0 \mathcal{S}} x$，即 $(x_d)_{d \in D} \to_{\mathcal{GI}_0 \mathcal{S}} x$。由假设知 $\{j \in J : x_j \notin \uparrow G\} \in \mathcal{I}_0$。因为 $\mathcal{I}_0$ 是非平凡理想，则 $\{j \in J : x_j \notin \uparrow G\} \neq D$。因此存在 $x_{d_0} = d_0$ 使得 $d_0 \in \uparrow G$，即 $G \ll x$。

**命题 6** 设 $L$ 是拟连续 domain，$G \in \mathcal{P}^w(L)$，$(x_j)_{j \in J} \subseteq L$ 是一个网，$\mathcal{I}$ 是 $J$ 的非平凡理想，$x \in L$。如果任意 $G \ll x$ 蕴含 $\{j \in J : x_j \notin \uparrow G\} \in \mathcal{I}$，则 $(x_j)_{j \in J} \to_{\mathcal{GIS}} x$。

**证明** 设 $G \in \mathcal{P}^w(L)$，$(x_j)_{j \in J} \subseteq L$ 是一个网，$\mathcal{I}$ 是 $J$ 的非平凡理想，$x \in L$ 且任意 $G \ll x$ 蕴含 $\{j \in J : x_j \notin \uparrow G\} \in \mathcal{I}$。由拟连续性知 $fin(x) = \{F \in \mathcal{P}^w(L) : F \ll x\}$ 是定向集族，且 $\uparrow x = \bigcap_{F \in \mathcal{F}} \uparrow F$。由假设知任意 $F \in fin(x)$，$\{j \in J : x_j \notin \uparrow F\} \in \mathcal{I}$，则由定义 4 知 $(x_j)_{j \in J} \to_{\mathcal{GIS}} x$。

**命题 7** 设 $L$ 是定向完备偏序集，则 $\tau_L^{g-I-\lim\text{-inf}}$ 是 $L$ 上的拓扑，称为广义理想下极限拓扑。

**证明** 这里只证明开集的交仍然是开集，其他显然。设 $U, V \in \tau_L^{g-I-\lim\text{-inf}}$，$(x_j)_{j \in J} \to_{\mathcal{GIS}} x$，$x \in U \cap V$，则

$$\{j \in J : x_j \notin U \cap V\} = \{j \in J : x_j \notin U\} \cup \{j \in J : x_j \notin V\} \in \mathcal{I},$$

则 $U \cap V \in \tau_L^{g-I-\lim\text{-inf}}$。

**命题 8** 设 $L$ 是定向完备偏序集，如果 $(x_j)_{j \in J} \to_{\mathcal{GIS}} x$，则 $\tau_L^{g-I-\lim\text{-inf}}$ 是使得 $(x_j)_{j \in J} \to_{\mathcal{I}} x$ 的最细的

拓扑。

**证明** 设$(x_j)_{j\in J}\to_\mathcal{I} x$关于拓扑$\tau$成立，下面证明$\tau\subseteq\tau_L^{g\text{-}I\text{-}\lim\text{-}\inf}$. 任意$x\in U\in\tau$，由$(x_j)_{j\in J}\to_\mathcal{I} x$知$\{j\in J:x_j\notin U\}\in\mathcal{I}$. 由$(x_j)_{j\in J}\to_{\mathcal{GIS}} x$，$x\in U\in\tau$和$\{j\in J:x_j\notin U\}\in\mathcal{I}$知$U\in\tau_L^{g\text{-}I\text{-}\lim\text{-}\inf}$.

设$L$是定向完备偏序集，$(x_j)_{j\in J}$是$L$中的一个网，$x\in L$. 如果存在定向集族$\mathcal{F}\in\mathcal{P}^w(X)$使得下列条件成立：(1)$\bigcap_{F\in\mathcal{F}}\uparrow F\subseteq\uparrow x$；(2)任意的$F\in\mathcal{F}$，存在$k\in J$使得当$j\geq k$时$x_j\in\uparrow F$. 则称网$(x_j)_{j\in J}$是广义下极限收敛到$x$，记为$(x_j)_{j\in J}\to_{\mathcal{GS}} x$ [15].

广义下极限收敛可以看成广义理想下极限收敛的特殊形式，因为$(x_j)_{j\in J}\to_{\mathcal{GS}} x$当且仅当$(x_j)_{j\in J}\to_{\mathcal{GI_0S}} x$. 首先如果$(x_j)_{j\in J}\to_{\mathcal{GS}} x$，则存在定向集族$\mathcal{F}\in\mathcal{P}^w(X)$使得$\bigcap_{F\in\mathcal{F}}\uparrow F\subseteq\uparrow x$且任意的$F\in\mathcal{F}$，存在$k\in J$使得当$j\geq k$时$x_j\in\uparrow F$，而
$$\{j\in J:x_j\notin\uparrow F\}=J\setminus\{j\in J:x_j\in\uparrow F\}\in\mathcal{I}_0$$
则$(x_j)_{j\in J}\to_{\mathcal{GI_0S}} x$；反之，如果$(x_j)_{j\in J}\to_{\mathcal{GI_0S}} x$，则存在定向集族$\mathcal{F}\in\mathcal{P}^w(X)$使得$\bigcap_{F\in\mathcal{F}}\uparrow F\subseteq\uparrow x$且任意的$F\in\mathcal{F}$，$\{j\in J:x_j\notin\uparrow F\}\in\mathcal{I}_0$，则存在$k\in J$使得$\{j\in J:x_j\notin\uparrow F\}\subseteq J\setminus M_k$，从而$M_k\subseteq\{j\in J:x_j\in\uparrow F\}$，即当$j\geq k$时$x_j\in\uparrow F$，$(x_j)_{j\in J}\to_{\mathcal{GS}} x$.

上述过程也说明对$J$的非平凡理想$\mathcal{I}_0$，如果$\{j\in J:x_j\notin\uparrow F\}\in\mathcal{I}_0$，则$(x_j)_{j\in J}$终在$\mathcal{I}_0$中. 对其它理想，上述结果不一定成立.

令
$$\tau_L^{g\text{-}\lim\text{-}\inf}=\{U\subseteq L:(x_j)_{j\in J}\to_{\mathcal{GS}} x,\ x\in U,\ 则(x_j)_{j\in J}\text{终在}U\text{中}\}.$$

易知$\tau_L^{g\text{-}\lim\text{-}\inf}$是$L$上的一个拓扑，称为广义下极限拓扑[15].

**命题 9** 设$L$是定向完备偏序集. 则$U\in\tau_L^{g\text{-}\lim\text{-}\inf}$当且仅当任意定向集族$\mathcal{F}\in\mathcal{P}^w(X)$满足$\bigcap_{F\in\mathcal{F}}\uparrow F\subseteq\uparrow x$，$x\in U$蕴含存在$F\in\mathcal{F}$使得$\uparrow F\subseteq U$.

**证明** 设$U\in\tau_L^{g\text{-}\lim\text{-}\inf}$，定向集族$\mathcal{F}\in\mathcal{P}^w(X)$满足$\bigcap_{F\in\mathcal{F}}\uparrow F\subseteq\uparrow x$，$x\in U$. 设任意的$F\in\mathcal{F}$，存在$x_F\in\uparrow F$使得$x_F\notin U$，显然网$(x_F)_{F\in\mathcal{F}}$广义下极限收敛到$x$. 由$\tau_L^{g\text{-}\lim\text{-}\inf}$的定义知网$(x_F)_{F\in\mathcal{F}}$关于广义下极限拓扑收敛，则$(x_j)_{j\in J}$终在$U$中，与$x_F\notin U$矛盾.

反之，如果任意定向集族$\mathcal{F}\in\mathcal{P}^w(X)$使得$\bigcap_{F\in\mathcal{F}}\uparrow F\subseteq\uparrow x$，$x\in U$蕴含存在$F\in\mathcal{F}$使得$\uparrow F\subseteq U$. 令$(x_j)_{j\in J}\to_{\mathcal{GS}} x$，则存在定向集族$\mathcal{F}\in\mathcal{P}^w(X)$使得$\bigcap_{F\in\mathcal{F}}\uparrow F\subseteq\uparrow x$，则由假设知存在$F\in\mathcal{F}$使得$\uparrow F\subseteq U$. 再由存在$k\in J$使得当$j\geq k$时$x_j\in\uparrow F$知$x_j\in U$，则$(x_j)_{j\in J}$终在$U$中，$U\in\tau_L^{g\text{-}\lim\text{-}\inf}$.

**命题 10** 设$L$是定向完备偏序集. 则$\tau_L^{g\text{-}\lim\text{-}\inf}=\tau_L^{g\text{-}I\text{-}\lim\text{-}\inf}$.

**证明** 首先证明$\tau_L^{g\text{-}I\text{-}\lim\text{-}\inf}\subseteq\tau_L^{g\text{-}\lim\text{-}\inf}$. 令$U\in\tau_L^{g\text{-}I\text{-}\lim\text{-}\inf}$，$x\in U$，网$(x_j)_{j\in J}$满足$(x_j)_{j\in J}\to_{\mathcal{GS}} x$，则$(x_j)_{j\in J}\to_{\mathcal{GI_0S}} x$，$\{j\in J:x_j\notin U\}\in\mathcal{I}_0$. 由$\mathcal{I}_0$的非平凡性知存在$x_k\in U$，当$j\geq k$时$x_j\in U$，$(x_j)_{j\in J}$终在$U$中，从而$U\in\tau_L^{g\text{-}\lim\text{-}\inf}$.

下面证明$\tau_L^{g\text{-}\lim\text{-}\inf}\subseteq\tau_L^{g\text{-}I\text{-}\lim\text{-}\inf}$. 设$U\in\tau_L^{g\text{-}\lim\text{-}\inf}$，$x\in U$，$x\leq y$，则$\{y\}$是定向集，且$\uparrow y=\bigcap_{\{y\}}\uparrow\{y\}$且$x\leq y$，即$(y)\to_{\mathcal{GS}} x$，则$(y)$终在$U$中，$y\in U$，$U=\uparrow U$. 如果$(x_j)_{j\in J}\to_{\mathcal{GIS}} x\in U$，则存在定向集族$\mathcal{F}\in\mathcal{P}^w(X)$使得$\bigcap_{F\in\mathcal{F}}\uparrow F\subseteq\uparrow x$且任意的$F\in\mathcal{F}$，$\{j\in J:x_j\notin\uparrow F\}\in\mathcal{I}$. 则由$\bigcap_{F\in\mathcal{F}}\uparrow F\subseteq\uparrow x\subseteq U$知，存在$F\in\mathcal{F}$使得$\uparrow F\subseteq U$，$\{j\in J:x_j\notin U\}\subseteq\{j\in J:x_j\notin\uparrow F\}\in\mathcal{I}$，即$\{j\in J:x_j\notin U\}\in\mathcal{I}$，$U\in\tau_L^{g\text{-}I\text{-}\lim\text{-}\inf}$.

上述命题说明了广义下极限拓扑与广义理想下极限拓扑一致. 下面说明在它们与偏序集上的Scott拓扑一致.

**定理 1** 设$L$是定向完备偏序集。则$\tau_L^{g\text{-}\lim\text{-}\inf}=\tau_L^{g\text{-}I\text{-}\lim\text{-}\inf}=\sigma(L)$.

**证明** 只需证明$\tau_L^{g\text{-}\lim\text{-}\inf}=\sigma(L)$. 设$U\in\sigma(L)$，如果$(x_j)_{j\in J}\to_{\mathcal{GS}} x\in U$，则存在定向集族

$\mathcal{F} \in \mathcal{P}^w(X)$ 使得 $\bigcap_{F \in \mathcal{F}} \uparrow F \subseteq \uparrow x$ 且任意的 $F \in \mathcal{F}$，存在 $k \in J$ 使得当 $j \geq k$ 时 $x_j \in \uparrow F$，则 $\bigcap_{F \in \mathcal{F}} \uparrow F \subseteq \uparrow x \subseteq U$，由 Rudin 引理知存在 $F \in \mathcal{F}$ 使得 $\uparrow F \subseteq U$，故 $(x_j)_{j \in J}$ 终在 $U$ 中，则 $U \in \tau_L^{g-\lim-\inf}$. 反之，如果 $U \in \tau_L^{g-\lim-\inf}$，$x \in U$，$x \leq y$，显然 $U$ 是上集. 设定向集 $D \subseteq L$ 满足 $\sup D \in U$，令 $\mathcal{F} = \{\{d\} : d \in D\}$，则 $\uparrow \sup D = \bigcap_{\{d\} \in \mathcal{F}} \uparrow d$ 且 $x \leq \sup D$，即 $\bigcap_{\{d\} \in \mathcal{F}} \uparrow d \subseteq \uparrow x \subseteq U$，则 $x \in U$ 蕴含存在 $\{d\} \in \mathcal{F}$ 使得 $\uparrow d \subseteq U$，则 $D \cap U \neq \varnothing$，$U \in \sigma(L)$.

文献[15]中命题 2.7 证明了类似结论，上述命题证明 $\sigma(L) \subseteq \tau_L^{g-\lim-\inf}$ 的过程与 $\tau_L^{g-\lim-\inf}$ 类似，$\tau_L^{g-\lim-\inf} \subseteq \sigma(L)$ 的证明过程本文采用了更直接的证明.

设 $L$ 是拟连续 domain，则 $\Uparrow F = \{x : F \ll x\}$ 是广义理想下极限拓扑 $\tau_L^{g-I-\lim-\inf}$ 的基.

**定义 5** 设 $L$ 是定向完备偏序集，$x \in L$，$\mathcal{I}$ 是 $J$ 的理想. 如果存在拓扑 $\tau$ 使得任意的网 $(x_j)_{j \in J} \subseteq L$，$(x_j)_{j \in J} \to_{\mathcal{GIS}} x$ 当且仅当 $(x_j)_{j \in J} \to_{\mathcal{I}} x$ 关于拓扑 $\tau$ 成立，则称广义理想下极限收敛是拓扑的.

**定理 2** 设 $L$ 是定向完备偏序集，则广义理想下极限收敛是拓扑的当且仅当 $L$ 是拟连续偏序集.

**证明** 令 $L$ 是拟连续偏序集. 设 $(x_j)_{j \in J} \to_{\mathcal{GIS}} x$，则显然 $(x_j)_{j \in J}$ 关于拓扑 $\tau_L^{g-I-\lim-\inf}$ 是理想收敛的. 设 $(x_j)_{j \in J} \to_{\mathcal{I}} x$，如果 $\mathcal{I}$ 是 $J$ 的平凡理想，则显然 $(x_j)_{j \in J} \to_{\mathcal{GIS}} x$. 设 $\mathcal{I}$ 是 $J$ 的非平凡理想，$(x_j)_{j \in J}$ 关于拓扑 $\tau_L^{g-I-\lim-\inf}$ 是理想收敛的. 因为 $L$ 是拟连续偏序集，则 $fin(x) = \{F \in \mathcal{P}^w(L) : F \ll x\}$ 是定向集族且 $\uparrow x = \bigcap_{F \in fin(x)} \uparrow F$. 任取 $F \in fin(x)$，$x \in \Uparrow F \in \tau_L^{g-I-\lim-\inf}$，因为 $(x_j)_{j \in J}$ 关于拓扑 $\tau_L^{g-I-\lim-\inf}$ 是理想收敛的，则 $\{j \in J : x_j \notin \Uparrow F\} \in \mathcal{I}$，则 $\{j \in J : x_j \notin \uparrow F\} \subseteq \{j \in J : x_j \notin \Uparrow F\} \in \mathcal{I}$，即 $(x_j)_{j \in J} \to_{\mathcal{GIS}} x$。

设广义理想下极限收敛是拓扑的，则存在拓扑 $\tau$ 使得任意的网 $(x_j)_{j \in J} \subseteq L$，$x \in L$，$(x_j)_{j \in J} \to_{\mathcal{GIS}} x$ 当且仅当 $(x_j)_{j \in J} \to_{\mathcal{I}} x$ 关于拓扑 $\tau$ 成立. 任意 $x \in L$，考察其在拓扑 $\tau$ 中的开邻域 $\mathcal{N}(x) = \{U \in \tau : x \in U\}$. 令 $M(x) = \{(U, y) : U \in \mathcal{N}(x), y \in U\}$，定义如下预序：$(U_1, y_1) \geq (U_2, y_2) \Leftrightarrow U_1 \subseteq U_2$，显然 $M(x)$ 是定向集. 设 $(x_{(U,y)})_{(U,y) \in M(x)} = y$，则 $(x_{(U,y)})_{(U,y) \in M(x)}$ 作为网满足 $(x_j)_{j \in J} \to_{\mathcal{I}_0} x$，即 $(x_j)_{j \in J} \to_{\mathcal{GI}_0 S} x$，则存在定向集族 $\mathcal{F} \subseteq \mathcal{P}^w(X)$ 使得 $\bigcap_{F \in \mathcal{F}} \uparrow F \subseteq \uparrow x$ 且任意的 $F \in \mathcal{F}$，有 $\{(U, y) \in M(x) : x_{(U,y)} \notin \uparrow F\} \in \mathcal{I}$. 由于 $\mathcal{I}_0$ 是 $J$ 的非平凡理想，则存在 $(U_F, y_F) \in M(x)$ 使得 $y_F = x_{(U_F, y_F)} \in \uparrow F$，则当 $x_{(U,y)} \geq x_{(U_F, y_F)}$ 时，$x_{(U,y)} \in \uparrow F$. 任取 $w \in U_F$，则 $(U_F, w) \geq (U_F, y_F)$，因此 $w = x_{(U_F, w)} \geq x_{(U_F, y_F)} \geq y_F$，即 $U_F \subseteq \uparrow F$ 对任意的 $F$ 成立.

设定向集 $D$ 满足 $x \leq \sup D$，因为 $\tau \subseteq \tau_L^{g-I-\lim-\inf}$，则 $U_F \in \tau_L^{g-I-\lim-\inf}$，即 $U_F$ 是 Scott 开集. 令 $\mathcal{F}_D = \{\{d\} : d \in D\}$，则存在 $\{d\} \in \mathcal{F}_D$ 使得 $\uparrow d \subseteq U_F$，故 $\uparrow d \subseteq \uparrow F$，$D \cap \uparrow F \neq \varnothing$，$F \ll x$. 则 $\mathcal{F} \subseteq fin(x)$ 且

$$\bigcap_{F \in fin(x)} \uparrow F \subseteq \bigcap_{F \in \mathcal{F}} \uparrow F = \uparrow a \subseteq \uparrow x,$$

显然 $\uparrow x \subseteq \bigcap_{F \in fin(x)} \uparrow F$，即 $\uparrow x = \bigcap_{F \in fin(x)} \uparrow F$. 则 $L$ 是拟连续偏序集.

## 3 广义理想终下界极限与 Lawson 拓扑

本节研究广义理想下极限收敛的特殊情况—广义理想终下界极限，并讨论其生成的拓扑与 Lawson 拓扑的关系. 最后给出了定向完备偏序集为连续 domain 的一个等价刻画.

设 $L$ 是定向完备偏序集，$(x_j)_{j \in J}$ 是 $L$ 中的一个网，$\mathcal{I}$ 是 $J$ 的理想，$x \in L$. 令

$$GI(x_j) = \{F \in \mathcal{P}^w(L) : \{j \in J : x_j \notin \uparrow F\} \in \mathcal{I}\}.$$

**定义 6** 设 $L$ 是定向完备偏序集，$(x_j)_{j \in J}$ 是 $L$ 中的一个网，$\mathcal{I}$ 是 $J$ 的理想，$x \in L$. 如果 (1) $(x_j)_{j \in J} \to_{\mathcal{GIS}} x$；(2) 任意 $\{j \in J : x_j \notin \uparrow F\} \in \mathcal{I}$，$x \in \uparrow F$，则称 $x$ 是网 $(x_j)_{j \in J}$ 广义理想终下界极限. 记为 $x = \underline{\lim}_{\mathcal{GI}}(x_j)_{j \in J}$.

设 $\mathcal{I}$ 是 $J$ 的理想，令

$$\tau_L^{g-I} = \{U \subseteq L : x = \underline{\lim}_{\mathcal{GI}}(x_j)_{j \in J}, x \in U，则 \{j \in J : x_j \notin U\} \in \mathcal{I}\}.$$

**命题 11** 设 $L$ 是定向完备偏序集，则 $\tau_L^{g-I}$ 是 $L$ 上的拓扑，称为广义理想终下界极限拓扑.

**命题 12** 设 $L$ 是定向完备偏序集，则 $\lambda(L) \subseteq \tau_L^{g-I}$.

**证明** 只需要证明 Scott 开集和下拓扑中的开集都是广义理想终下界极限拓扑中的开集.

设 $x = \underline{\lim}_{\mathcal{GI}}(x_j)_{j\in J}$ 且 $x \in U \in \sigma(L)$. 由 $x = \underline{\lim}_{\mathcal{GI}}(x_j)_{j\in J}$ 的定义知，存在定向集族 $\mathcal{F}$ 使得 $\mathcal{F} \subseteq GI(x_j)$ 使得 $\bigcap_{F\in\mathcal{F}} \uparrow F \subseteq \uparrow x$，则 $\bigcap_{F\in\mathcal{F}} \uparrow F \subseteq U$. 由 Rudin 引理及其推论知存在 $F \in \mathcal{F}$ 满足 $\uparrow F \subseteq U$. 对任意 $F \in \mathcal{F}$，$\{j \in J : x_j \notin \uparrow F\} \in \mathcal{I}$. 而 $\{j \in J : x_j \notin U\} \subseteq \{j \in J : x_j \notin \uparrow F\} \in \mathcal{I}$，$U \in \tau^{g-I}$.

要证明 $L \setminus \uparrow x$ 是 $\tau_L^{g-I}$ 中开集，只需证明 $\uparrow x$ 为 $\tau_L^{g-I}$ 中闭集. 设网 $(x_j)_{j\in J}$ 满足 $(x_j)_{j\in J} \subseteq \uparrow x$，则 $\{j \in J : x_j \notin \uparrow x\} = \varnothing \in \mathcal{I}$. 假设 $a = \underline{\lim}_{\mathcal{GI}}(x_j)_{j\in J}$，而 $\{j \in J : x_j \notin \uparrow x\} = \varnothing \in \mathcal{I}$，则 $x \leq a$. 由此可知 $\uparrow x$ 为 $\tau_L^{g-I}$ 中闭集，$L \setminus \uparrow x$ 是 $\tau_L^{g-I}$ 中开集.

综上所述可知 $\lambda(L) \subseteq \tau_L^{g-I}$.

**定理 3** 设 $L$ 是拟连续 domain，$\mathcal{I}$ 是 $J$ 的非平凡理想. 则 $x = \underline{\lim}_{\mathcal{GI}}(x_j)_{j\in J} \Leftrightarrow (x_j)_{j\in J}$ 在 Lawson 拓扑中满足 $(x_j)_{j\in J} \to_{\mathcal{I}} x$.

**证明** 一方面，如果 $x = \underline{\lim}_{\mathcal{GI}}(x_j)_{j\in J}$，则 $(x_j)_{j\in J}$ 在 Lawson 拓扑中满足 $(x_j)_{j\in J} \to_{\mathcal{I}} x$.

另一方面，设 $(x_j)_{j\in J}$ 在 Lawson 拓扑中满足 $(x_j)_{j\in J} \to_{\mathcal{I}} x$，因为 $L$ 是拟连续 domain，则定向集族 $fin(x) = \{F \in \mathcal{P}^w(L) : F \ll x\}$，$\bigcap_{F \in fin(x)} \uparrow F = \uparrow x$. 任意的 $F \ll x$，由 $L$ 的拟连续性可知 $\Uparrow F$ 为 Scott 开集，即 Lawson 开集. 又 $(x_j)_{j\in J}$ 关于 Lawson 拓扑收敛到 $x$，则 $(x_j)_{j\in J}$ 终在 $\Uparrow F$ 中. 因为 $x \in \Uparrow F \subseteq \uparrow F$ 且 $(x_j)_{j\in J}$ 终在 $\Uparrow F$ 中，则 $\{j \in J : x_j \notin \Uparrow F\} \in \mathcal{I}$，即 $\{j \in J : x_j \notin \uparrow F\} \subseteq \{j \in J : x_j \notin \Uparrow F\} \in \mathcal{I}$.

任意 $\{j \in J : x_j \notin \uparrow F\} \in \mathcal{I}$，如果 $x \notin \uparrow F$，则 $x \in L \setminus \uparrow F$，而 $L \setminus \uparrow F$ 是 Lawson 开集且 $(x_j)_{j\in J}$ 在 Lawson 拓扑中满足 $(x_j)_{j\in J} \to_{\mathcal{I}} x$，则 $\{j \in J : x_j \notin L \setminus \uparrow F\} \in \mathcal{I}$，而 $\{j \in J : x_j \notin \uparrow F\} \in \mathcal{I}$，即 $J = \{j \in J : x_j \notin L \setminus \uparrow F\} \cup \{j \in J : x_j \notin \uparrow F\} \in \mathcal{I}$，与 $\mathcal{I}$ 是 $J$ 的非平凡理想矛盾，因此 $x \in \uparrow F$. 综上可知 $x = \underline{\lim}_{\mathcal{GI}}(x_j)_{j\in J}$.

**定理 4** 设 $L$ 是交连续的定向完备偏序集，$\mathcal{I}$ 是 $J$ 的非平凡理想. 如果 $x = \underline{\lim}_{\mathcal{GI}}(x_j)_{j\in J} \Leftrightarrow (x_j)_{j\in J}$ 在 Lawson 拓扑中满足 $(x_j)_{j\in J} \to_{\mathcal{I}} x$ 且任意 $\{j \in J : x_j \notin \uparrow F\} \in \mathcal{I}$ 蕴含 $(x_j)_{j\in J}$ 终在 $\uparrow F$，则 $L$ 是连续的.

**证明** 设 $x = \underline{\lim}_{\mathcal{GI}}(x_j)_{j\in J} \Leftrightarrow (x_j)_{j\in J}$ 在 Lawson 拓扑中满足 $(x_j)_{j\in J} \equiv_{\mathcal{I}} x$. 任意 $x \in L$，考察其在 Lawson 拓扑中的开邻域 $\mathcal{N}(x) = \{U \in \tau : x \in U\}$. 定义 $I = \{(U, n, a) \in \mathcal{N}(x) \times N \times L : a \in U\}$，在 $I$ 上定义序：$(U, n, a) \leq (V, m, b) \Leftrightarrow V \subseteq U, V \neq U$ 或 $V = U, n \leq m$. 对任意 $i = (U, n, a) \in I$，令 $x_i = a$. 这样定义的网 $(x_i)_{i \in I}$ 关于 Lawson 拓扑收敛到 $x$，所以 $x = \underline{\lim}_{\mathcal{GI}}(x_i)_{i \in I}$. 则存在定向集族 $\mathcal{F}$ 使得 $\bigcap_{F \in \mathcal{F}} \uparrow F \subseteq \uparrow x$. 令 $F \in \mathcal{F}$，则 $\{j \in J : x_j \notin \uparrow F\} \in \mathcal{I}$. 存在 $i = (U, n, a) \in I$，由于对 $(V, m, b) = j \geq i$ 有 $x_j \in \uparrow F$. 特别地，对所有 $b \in U$，$j = (U, m+1, b) \geq (U, m, a)$ 都满足 $x_j = b \in \uparrow F$，所以 $U \subseteq \uparrow F$. 由于 $U$ 为 Lawson 开集，$x \in U$ 则存在 Scott 开集 $V$ 及有限集 $W$ 满足 $x \in V \setminus \uparrow W \subseteq U \subseteq \uparrow F$，所以 $x \in \uparrow (V \setminus \uparrow W)$. 由于 $V \setminus \uparrow W$ 是 Lawson 开集且 $L$ 是交连续的，则 $\uparrow (V \setminus \uparrow W)$ 是 Scott 开集，从而 $x \in (\uparrow F)^\circ_{\sigma(L)}$. 对任意定向集 $D$ 满足 $x \leq \sup D$，有 $\sup D \in (\uparrow F)^\circ_{\sigma(L)}$. 于是存在 $d \in D$ 满足 $d \in (\uparrow F)^\circ_{\sigma(L)} \subseteq \uparrow F$，即 $F \ll x$. 注意到 $\uparrow a = \bigcap_{F \in \mathcal{F}} \uparrow F$，$x \leq a$，即 $\bigcap_{F \in \mathcal{F}} \uparrow F \subseteq \uparrow x$. 从而 $L$ 是拟连续的，所以 $L$ 是连续的.